\documentclass[11pt]{article}

\usepackage{epic,eepic}
\usepackage{amscd,amsfonts,amsmath,amsxtra,amssymb}
\usepackage[all]{xy}
\input xypic
\newfont{\gothic}{eufm10 scaled 1100}
\newfont{\mediumgothic}{eufm10 scaled 1000}
\newfont{\smallgothic}{eufm10 scaled 900}
\newfont{\verysmallgothic}{eufm10 scaled 700}

\def \PP {{\mathbb P}}

\def \ZZ {{\mathbb Z}}

\newcommand{\mgamma}{{m/\Gamma}}
\newcommand{\oka}{{\cal O}}
\newcommand{\ms}{{{M}}}     					
\newcommand{\msbar}{\overline{{M}}}				
\newcommand{\mgstack}[1]{{\overline{\cal M}_{{#1}} }}		
\newcommand{\mgstackopen}[1]{{{\cal M}_{{#1}}}}			
\newcommand{\Hgn}{\overline{H}_{g,n,m}}
\newcommand{\Hgnopen}{{H}_{g,n,m}}
\newcommand{\Hgmp}{\overline{H}_{g,n,\mgamma}}			

\renewcommand{\epsilon}{\varepsilon}
\renewcommand{\rho}{\varrho}
\renewcommand{\theta}{\vartheta}
\newcommand{\ob}{\text{\rm Ob}\,} 
\newcommand{\simto}{\stackrel{\sim}{\longrightarrow}}

\newcommand{\Hilb}{\text{\rm Hilb}\,}

\newcommand{\Aut}{\text{\rm Aut}  }
\newcommand{\Stab}{\text{\rm Stab}}
\newcommand{\id}{\text{\rm id}}
\newcommand{\pr}{\text{\rm pr}}

\newcommand{\PGL}{\text{\rm PGL}}
\newcommand{\Spec}{\text{\rm Spec\,}}

\newcommand{\ebew}{\hfill$\Box$ \par}

\newcommand{\proof}{{\em \underline{Proof.}\quad}}
\newcommand{\proofof}[1]{{\em \underline{Proof of #1.}\quad}}
\newtheorem{thm}{Theorem}[section]
\newtheorem{defi}[thm]{Definition}
\newtheorem{prop}[thm]{Proposition}

\newtheorem{rk}[thm]{Remark}

\newtheorem{lemma}[thm]{Lemma}

\newtheorem{abschnitt}[thm]{}

\begin{document}

\title{Moduli stacks\\ of permutation classes of pointed stable curves}

\author{J\"org Zintl}
\date{11.11.2006}
\maketitle

\begin{abstract}
The notion of $\mgamma$-pointed stable curves is introduced. It should be viewed as a generalization of the notion of $m$-pointed stable curves of a given genus, where the labels of the marked points are only determined up to the action of a group of permutations $\Gamma$. The classical moduli spaces and moduli stacks are generalized to this wider setting. Finally, an explicit construction of the new moduli stack of  $\mgamma$-pointed stable curves as a quotient stack is given. 
\end{abstract}

\section{Introduction}
By a permutation class of a pointed stable curve over a closed point $\Spec(k)$ we mean a stable curve with marked points, where the labels of the marked points are only determined up to certain permutations.

More formally, let $m\in{\mathbb N}$ be a natural number, and let $\Gamma$ be a subgroup of the group $\Sigma_m$ of permutations on $m$ elements. 
An $m$-pointed stable curve of genus $g$ over a closed point is given by a tuple $(f:C\rightarrow\Spec(k);P_1,\ldots,P_m)$, where $P_1,\ldots,P_m$ define $m$ marked smooth points on an algebraic curve $C$, compare definition \ref{mpointed}. A permutation class of pointed stable curves with respect to $\Gamma$  is simply an equivalence class of such tuples, where  $(f:C\rightarrow\Spec(k),P_1,\ldots,P_m)$ and $(f':C'\rightarrow\Spec(k);P_1',\ldots,P_m')$ are called equivalent if $C=C'$  and if  there exists a permutation $\gamma\in\Gamma$, such that $P_i = P_{\gamma(i)}'$ for all $i=1,\ldots,m$. Below, such a class will be called an $\mgamma$-pointed stable curve over $\Spec(k)$, see definition \ref{mgdef} and lemma \ref{onetoone04}.  

Over an arbitrary base scheme $S$ the situation is more subtle. It turns out that an $\mgamma$-pointed stable curve cannot simply be defined as an equivalence class of $m$-pointed stable curves with respect to the analogous equivalence relation.  However, the definition of $\mgamma$-pointed stable curves is made in such a way that locally, with respect to the \'etale topology,    
an $\mgamma$-pointed stable curve $f:C\rightarrow S$ over an arbitrary scheme $S$ looks like an equivalence class  of $m$-pointed stable curves over $S$. The  sections of marked points, which exist on an \'etale cover, are only determined up to permutations in $\Gamma$, and they need to satisfy certain compatibility conditions, but in general they will not descent to global sections of the original curve. 

The classical $m$-pointed stable curves are special cases of $\mgamma$-pointed stable curves, where $\Gamma$ is the trivial group. 

More generally, stable curves with $m$ distinguished points, where the points are only ordered up to certain permutations, occur on a number of occasions. Consider for example a stable curve $f: C\rightarrow \Spec(k)$ of some genus $g$, which consists of more than one irreducible component. An irreducible component of $C$ is a curve $f': C'\rightarrow \Spec(k)$ of some genus $g' \le g$, but it no longer needs to  be a stable curve. However, the curve $f': C'\rightarrow \Spec(k)$ together with the unordered set $\{P_1,\ldots,P_m\}$, consisting of those  closed points where $C'$ intersects the closure of its comlement in $C$, is stable when considered as an $\mgamma$-pointed stable curve of genus $g'$, where $\Gamma$ is the group of all permutations on $m$ elements.   
Note that even if an ordering of the nodes of $C$ is fixed a priory, there is in general no distinguished ordering of the intersection points $\{P_1,\ldots,P_m\}$. Because of automorphisms of $C$, an ordering will be determined only up to certain permutations. 

Another example, where permutation classes of pointed stable curves show up naturally, is the moduli space  $\overline{M}_{g}$ of Deligne-Mumford stable curves of genus $g$ itself.  There is a canonical stratification $\overline{M}_{g} = \bigcup_i {M}_{g}^{(i)}$, where ${M}_{g}^{(i)}$ denotes the locus of curves with exactly $i$ nodes, for $i=0,\ldots,3g-3$. The irreducible components of the subschemes of the stratification can be related to moduli spaces of $\mgamma$-pointed stable curves of genus $g'$,  where the number $m$, the group $\Gamma$ and the genus $g'$ varies. The case $i = 3g-4$ has been discussed in detail in \cite{Zi}, with special focus on the relation of the corresponding moduli stacks. 

Due to the fact that $\mgamma$-pointed stable curves over closed points can be described as equivalence classes of $m$-pointed stable curves, it is possible to construct coarse and fine moduli spaces for them  as quotients of moduli spaces of $m$-pointed stable curves. More importantly, because of the \'etale nature of $\mgamma$-pointed stable curves, we can define the corresponding moduli stacks in the sense of \cite{DM}. Finally, generalizing an idea of Edidin \cite{Ed} for Deligne-Mumford stable curves, we describe explicitely how to construct the new moduli stack of  $\mgamma$-pointed stable curves as a quotient stack.

\section{Preliminaries on pointed stable curves}
\begin{abschnitt}\em
We want to collect some basics about pointed stable curves first, so that we can refer to them later on. We also generalize a theorem of Edidin \cite{Ed} on the moduli stack of Deligne-Mumford stable curves to the case of pointed stable curves.\\

Let us first recall the definition of pointed stable curves as it is given in the paper of Knudsen  \cite[def. 1.1]{Kn}.
\end{abschnitt}

\begin{defi}\em\label{mpointed}
An {\em $m$-pointed stable curve of genus $g$} is a flat and proper
morphism $f : C \rightarrow S$ of schemes, together with $m$ sections
$\sigma_i : S \rightarrow C$, for $i=1,\ldots,m$, such that for all closed points $s\in
S$ holds
\begin{itemize}
\item[$(i)$] the fibre $C_s$ is a reduced connected algebraic curve with at most
ordinary double points as singularities; 

\item[$(ii)$] the arithmetic genus of $C_s$ is $\dim
H^1(C_s,\oka_{{C_s}}) = g$;

\item[$(iii)$] for $1\le i \le m$, the point $\sigma_i(s)$ is a smooth
point of $C_s$;

\item[$(iv)$] for all $ 1\le i,j \le m$ holds $\sigma_i(s) \neq
\sigma_j(s)$ if $i \neq j$;

\item[$(v)$] the number of points where a nonsingular rational
component $C_s'$ of $C_s$ meets the rest of $C_s$ plus the number of
points $\sigma_i(s)$ which lie on $C_s'$ is at least $3$.
\end{itemize}
\end{defi}
 
\begin{rk}\em
Condition $(v)$ guarantees that the group $\Aut(C_s)$ of automorphisms of a fibre $C_s$ is finite. Because of condition $(iv)$, the permutation group $\Sigma_m$\label{n3} acts faithfully on the set of sections $\{\sigma_1,\ldots,\sigma_m\}$, as well as on the set of marked points $\{\sigma_1(s),\ldots,\sigma_m(s)\}$ of each fibre $C_s$. 
\end{rk} 

\begin{abschnitt}\em
Let $f: C\rightarrow S$ be an $m$-pointed  stable curve of genus $g$ such that  $2g-2+m > 0$. The sections $\sigma_1,\ldots, \sigma_{m}$ determine effective Cartier divisors $S_1,\ldots ,S_{m}$ on $C$. 

There is a  canonical invertible 
dualizing sheaf on $C$, which is  denoted by $\omega_{C/S}$. In \cite[cor. 1.9, cor 1.11]{Kn} it is  shown that for
$n\ge 3$ the sheaf $\left(\omega_{C/S}(S_1+\ldots+S_{m})\right)^{\otimes n}$ is
relatively very ample,   and furthermore that $f_\ast\left(\left(\omega_{C/S}(S_1+\ldots+S_{m})\right)^{\otimes n}\right)$ is locally free of
rank $(2g-2+m)n -g+1$. 

Consider the Hilbert scheme $\Hilb_{{\PP^N}}^{{P_{g,n,m}}}$ of curves $ C \rightarrow \Spec(k)$ embedded in $\PP^N$ with Hilbert polynomial $P_{g,n,m} := (2g-2+m)nt-g+1$, where  $N :=(2g-2+m)n -g$. 

The Hilbert scheme of a simple point in $\PP^N$ is of course $\PP^N$ itself. The incidence condition of $m$ points lying on a curve $C$ defines a closed subscheme 
\[ I \subset 
\Hilb_{{\PP^N}}^{{P_{g,n,m}}} \times (\PP^N)^m .\]
There is an open subset $U \subset I$ parametrizing curves, where the $m$ points are pairwise different, and smooth points of $C$. Finally, there is a closed subscheme $ \overline{H}_{g,n,m} $ of $U$, which represents such embedded $m$-pointed stable curves $C \rightarrow \Spec(k)$ of genus $g$, where the embedding is determined by an isomorphism 
\[  \left(\omega_{C/\Spec(k)}(Q_1+\ldots+Q_{m})\right)^{\otimes n} \cong \oka_{{\PP^N}}(1) | C. \]
Here, $Q_1,\ldots,Q_m$ denote the marked points on $C$. 
The scheme constructed in this way 
is in fact a quasi-projective subscheme
\[ \overline{H}_{g,n,m} \subset \Hilb_{{\PP^N}}^{{P_{g,n,m}}} \times (\PP^N)^m \label{n3a}\]
of the Hilbert scheme. There is also a subscheme ${H}_{g,n,m }\subset\overline{H}_{g,n,m}$ corresponding to non-singular stable curves.  Compare \cite[section 2.3]{FP} for the more general situation of moduli of stable maps.  

There is a one-to-one correspondence  
between morphisms $\theta: S \rightarrow \Hgn$ and $m$-pointed stable curves $f :
C\rightarrow S$ together with global trivializations of Grothendieck's associated projective bundle 
\[ \PP f_\ast\left(\left(\omega_{C/S}(S_1+\ldots+S_{m})\right)^{\otimes n}\right) \: 
\simto \: \PP^N\times S. \]
In particular, the curve $f: C\rightarrow S$ can be considered as embedded into $\PP^N\times S$, such that the diagram 
\[\diagram
C \rrto|<\hole|<<\ahook \drto_f &&\PP^N \dlto^{\pr_2}\times S \\
&S 
\enddiagram \]
commutes. 

An element of the group $\PGL(N+1)$ acts on $\PP^N$, changing embeddings of $C$ by an isomorphism. Hence there is a natural action of $\PGL(N+1)$  on $\Hilb_{{\PP^N}}^{{P_{g,n,m}}}$, and thus on $\Hilb_{{\PP^N}}^{{P_{g,n,m}}} \times (\PP^N)^m$, which restricts to an action on 
$\Hgn$. Furthermore, the symmetric group $\Sigma_m$ acts on $(\PP^N)^m$ by permutation of the coordinates, which is equivalent to the permutation of the $m$ sections of marked points. Hence $\Sigma_m$ acts freely on $\Hgn\times (\PP^N)^m$. The combined action of  $\PGL(N+1) \times \Sigma_m$ on $\Hgn$ shall be written as an  action from the right. 
\end{abschnitt}

\begin{abschnitt}\em
Let $f: C_0 \rightarrow \Spec(k)$ be an $m$-pointed  stable curve over an algebraically closed field $k$, so in
particular an algebraic curve of genus $g$ with nodes as its only
singularities,  and $m$ marked points $P_1, \ldots , P_{m}$.  If we
fix one (arbitrary)  embedding of $C_0$ into $\PP^N$ via an isomorphism 
\[ \PP f_\ast(\omega_{C_0/k}(P_1+\ldots+P_m))^{\otimes n} \cong \PP^N ,\]
then this
distinguishes a $k$-valued point 
$[C_0]\in \Hgn$\label{n3b}. For the group of automorphisms of $C_0$ respecting the marked points, i.e. automorphisms which map each marked point to itself, there is a
natural isomorphism
\[ \Aut(C_0) \cong \Stab_{\PGL(N+1)}([C_0]),\]
that is,   an isomorphism with the subgroup of those  elements in $\PGL(N+1)$, which stabilize  the
point $[C_0]$. Recall that $\Sigma_m$ acts on $\Hgn$ without fixed points, so the stabilizer of $[C_0]$ in $\PGL(N+1)$ can be identified with the stabilizer of $[C_0]$ in $\PGL(N+1)\times \Sigma_m$. 

In general, elements $\gamma \in \PGL(N+1)$ are in
one-to-one correspondence with isomorphisms
\[ \gamma : \: \: C_0 \rightarrow C_\gamma\]
of embedded $m$-pointed stable curves, where 
 $C_\gamma$ is the curve
represented by the  point $[C_\gamma] := [C_0] \cdot \gamma$ in
$\Hgn$.
\end{abschnitt}

\begin{abschnitt}\em
Generalizing the construction of  $\msbar_{g}$ as described by Gieseker, the moduli space $\msbar_{g,m}$\label{n4}
of $m$-pointed stable
curves of genus $g$ is constructed as the GIT-quotient of the action
of $\PGL(N+1)$ on $\Hgn$, see \cite{Gs}, \cite{DM} and \cite{HM} for
the case of $m=0$, and compare \cite[Remark 2.4]{FP} for the more general case of stable maps. Note that the notion of stability of embedded $m$-pointed stable curves is compatible with the notion of stability of the corresponding points in $\Hgn$  with respect to the group action. 
In particular, there is a canonical quotient morphism
\[ \Hgn \rightarrow \msbar_{g,m} .\]
The construction of the moduli space is independent of the choice of $n$, provided it is large enough, and
hence independent of $N$.
For a stable algebraic curve $C_0 \rightarrow \Spec(k)$, the fibre over the point
in $\msbar_{g,m}$ representing it is isomorphic to the quotient $ \Aut(C_0) \setminus
\PGL(N+1)$.

The scheme $\msbar_{g,m}$ is at the same time a moduli space for the stack
$\mgstack{g,m}$\label{n5} of $m$-pointed stable curves of genus $g$.
\end{abschnitt}

\begin{defi}\em
The {\em moduli stack $\, \mgstack{g,m}$ of $\,m$-pointed stable curves of genus $g$} is the stack defined as the category fibred in groupoids over  the category of schemes, where for a scheme $S$ the objects in the fibre category $\mgstack{g,m}(S)$ are the $m$-pointed stable curves of genus $g$ over $S$.

 Morphisms in $\mgstack{g,m}$ are given as follows. Let $f: C \rightarrow S$ and $f': C' \rightarrow S'$ be $m$-pointed stable curves of genus $g$, i.e. objects of $\mgstack{g,m}(S)$ and $\mgstack{g,m}(S')$, respectively, with sections $\sigma_i:S \rightarrow C$ and $\sigma_i':S' \rightarrow C'$ for $i=1,\ldots,m$. A morphism from $f: C \rightarrow S$ to $f': C' \rightarrow S'$ in $\mgstack{g,m}$ is given by a pair of morphisms of schemes
\[ (g: S' \rightarrow S, \; \overline{g}: C' \rightarrow C),\]
such that the diagram
\[ \diagram
C' \rrto^{\overline{g}} \dto_{f'}&& C \dto^{f}\\
S' \rrto_g& &S
\enddiagram \]
is Cartesian, and for all $i=1,\ldots,m$ holds
\[ \overline{g}\circ \sigma_i' = \sigma_i \circ g .\]
\end{defi}

\begin{rk}\em
It was shown by Knudsen \cite{Kn} that the moduli stack $\mgstack{g,m}$ is a smooth, irreducible Deligne-Mumford stack, which is  proper over $\Spec(\ZZ)$. Compare also \cite{DM} and \cite{Ed}. 
\end{rk}

The following proposition was proven for the case $m=0$ in \cite{Ed}, but also holds true in the case of arbitrary $m$. Here, as usual, square brackets are used to denote quotient stacks. 

\begin{prop}\label{28}
There are  isomorphisms of stacks
\[ \mgstack{g,m} \cong \left[ \, \Hgn/\PGL(N+1) \right], \]
and
\[ \mgstackopen{g,m} \cong \left[ \, \Hgnopen/\PGL(N+1) \right]. \]
\end{prop}

\proof
We will only give the proof in the case of the closed moduli stack, the second 
 case being  completely analogous. 
Recall that for any given scheme $S$ an object of  the fibre category 
$\left[ \, \Hgn/\PGL(N+1) \right](S)$ is a 
triple $(E,p,\phi)$, where $p: E\rightarrow S$ is a principal
$\PGL(N+1)$-bundle, and $\phi: E\rightarrow \Hgn$ is a
$\PGL(N+1)$-equivariant morphism. 

$(i)$
An isomorphism of stacks from $\mgstack{g,m}$ to $ [\Hgn/\PGL(N+1) ]$ is
constructed as follows. Let an $m$-pointed  stable curve $f : C
\rightarrow S \in 
\ob(\mgstack{g,m})$ be given, 
with sections $\sigma_1,\ldots,\sigma_{m}: S \rightarrow C$ 
defining divisors $S_1,\ldots,S_{m}$ on $C$.   We denote  by $p :
E\rightarrow S$ the 
principal $\PGL(N+1)$-bundle associated to the projective bundle 
$\PP f_\ast(\omega_{C/S}(S_1+\ldots+S_{m}))^{\otimes
n}$ over $S$.  Consider the
Cartesian diagram
\[\diagram
C\times_S E \rto^{\overline{f}} \dto_{\overline{p}} & E \dto^p\\
C \rto_f & S.
\enddiagram \]
The pullback of $\PP f_\ast(\omega_{C/S}(S_1+\ldots+S_{m}))^{\otimes
n} $ to $E$ has a natural trivialization as 
a  projective bundle. Denote by $\tilde{S}_1,\ldots,\tilde{S}_{m}$ the divisors on $C\times_S E$ defined by the sections $\tilde{\sigma}_i : E \rightarrow C\times_S E$, with $\tilde{\sigma}_i(e) := [\sigma_i\circ p(e),e]$ for $e\in E$ and $i=1,\ldots,m$. Then by the universal property of the pullback there is a unique isomorphism  
\[ p^\ast  \PP f_\ast\left(\omega_{C/S}(S_1+\ldots+S_{m})\right)^{\otimes
n}  \cong 
\PP \overline{f}_\ast\left(\omega_{C\times_S E/E}   (\tilde{S}_1+\ldots+\tilde{S}_{m})\right)^{\otimes
n}  \]
of projective bundles over $E$. 
So we have an $m$-pointed  stable curve $ \overline{f}: C\times_S E
\rightarrow E$, together with a natural trivialization of the projective bundle $\PP
\overline{f}_\ast(\omega_{C\times_S E/E} (\tilde{S}_1+\ldots+\tilde{S}_{m}))^{\otimes
n}$. This is equivalent to specifying a morphism 
\[\phi : \: \: E \rightarrow \Hgn \]
by the universal property of the Hilbert scheme. By construction, the morphism  $\phi$
is $\PGL(N+1)$-equivariant, so we obtain an object $(E,p,\phi)\in
[\Hgn/\PGL(N+1) ](S)$.

To define a functor from $ \mgstack{g,m}$ to $ [\Hgn/\PGL(N+1) ]$ we need also to consider morphisms $(g:S'\rightarrow S, \overline{g}:C'\rightarrow C)$ between  stable curves $f: C\rightarrow S$ and $f':C'\rightarrow S'$. Because for the respective sections $ \sigma_1,\ldots,\sigma_m: S\rightarrow C$ and $\sigma_1',\ldots,\sigma_m':S'\rightarrow C'$ holds $\sigma_i\circ g = \overline{g}\circ \sigma_i'$ for all $i=1,\ldots,m$ by definition,  there is an induced morphism
\[ \PP f_\ast'(\omega_{C'/S'}(S_1'+\ldots+S_m'))^{\otimes n} \rightarrow  
\PP f_\ast(\omega_{C/S}(S_1+\ldots+S_m))^{\otimes n} ,\]
where $S_1',\ldots,S_m'$ denote the divisors of the marked points on $C'$. Hence there is an induced morphism of the associated principal $\PGL(N+1)$-bundles $\tilde{g} : E' \rightarrow E$, which fits into a commutative diagram
\[ \diagram
&& \Hgn\\
E' \rto_{\tilde{g}} \dto_{p'} \urrto^{\phi'} & E \dto^{p}\urto_{\phi}\\
S' \rto_g & S.
\enddiagram \]
One easily verifies that this assignment is functorial. In this way we obtain a functor from the fibred category $\mgstack{g,m}$ to the quotient fibred category $ [\Hgn/\PGL(N+1) ]$ over the category of schemes, and thus a morphism of stacks. 

$(ii)$ Conversely, consider a triple $(E,p,\phi)\in
[\Hgn/\PGL(N+1) ](S)$. The morphism $\phi : E \rightarrow \Hgn $
determines an $m$-pointed  stable curve $f' : C'\rightarrow E$ of genus $g$,
together with a trivialization $\PP f'_\ast(\omega_{C'/E}(S_1'+\ldots+S_{m}'))^{\otimes
n} \cong \PP^N \times E$, where $S_1',\ldots,S_{m}'$ denote the divisors on $C'$ determined by the $m$ sections $\sigma_1',\ldots,\sigma_{m}' : E \rightarrow C'$. 
The group  $\PGL(N+1)$ acts  diagonally on  $\PP^N \times E$. By the $\PGL(N+1)$-equivariance of $\phi$, we have for all $e\in E$ and all $\gamma \in \PGL(N+1)$ the identity of fibres $C'_e \cdot \gamma = C'_{\gamma(e)}$ as embedded $m$-pointed curves. Therefore there is an induced action of $\PGL(N+1)$ on the embedded curve $C'\subset \PP^N \times E$, which respects the sections $\sigma_1',\ldots,\sigma_{m}'$. 
Taking quotients we obtain 
\[ C := C' / \PGL(N+1) \longrightarrow E/ \PGL(N+1) \cong S, \]
which defines an $m$-pointed  stable curve $C \rightarrow S\in \mgstack{g,m}(S)$.

The same construction assigns to a morphism in $[\Hgn/\PGL(N+1) ]$ a morphism between $m$-pointed stable curves. This defines a functor between fibred categories, and thus
a morphism of stacks from  $[\Hgn/\PGL(N+1) ]$ to  $\mgstack{g,m}$. 

$(iii)$ It remains to show that the two functors defined above form an equivalence of categories. In fact, the composition
\[\mgstack{g,m} \rightarrow [\Hgn/\PGL(N+1) ]\rightarrow\mgstack{g,m} \]
is just the identity on $\mgstack{g,m}$.

Conversely, let $(E,p,\phi) \in [\Hgn/\PGL(N+1) ](S)$ for some scheme $S$. This defines an $m$-pointed  stable curve $f:C\rightarrow S$ in $\mgstack{g,m}(S)$, which in turn defines an object $(E',p',\phi') \in [\Hgn/\PGL(N+1) ](S)$. We claim that the triples $(E,p,\phi)$ and $(E',p',\phi')$ are isomorphic. 

By construction we have a commutative diagram
\[ \diagram
C' \rto^{f'} \dto_{\overline{p}} & E\dto^p\\
C := C'/\PGL(N+1) \rto_f & E/\PGL(N+1) = S, 
\enddiagram \]
where $f':C'\rightarrow E$ is the $m$-pointed stable curve over $E$ defined by the morphism $\phi:E\rightarrow \Hgn$, and the vertical arrows represent the natural quotient morphisms.  Hence there is an isomorphism of projective bundles
\[ p^\ast\PP f_\ast(\omega_{C/S}(S_1+\ldots +S_m))^{\otimes n} \cong 
\PP f'_\ast(\omega_{C'/E}(S'_1+\ldots +S'_m))^{\otimes n}, \]
together with a trivialization of the bundle $
\PP f'_\ast(\omega_{C'/E}(S'_1+\ldots +S'_m))^{\otimes n}$ determined by $\phi: E \rightarrow \Hgn$. 
Thus the pullback of the projective bundle $\PP f_\ast(\omega_{C/S}(S_1+\ldots+S_m))^{\otimes n}$ to the principal $\PGL(N+1)$-bundle $p: E \rightarrow S$ is trivial, and hence $E$ must be isomorphic to the principal $\PGL(N+1)$-bundle associated to $\PP f_\ast(\omega_{C/S}(S_1+\ldots+S_m))^{\otimes n}$.

It is straightforward to give the corresponding arguments for the 
morphisms in both categories to show that  the two constructions are inverse to each other, up to isomorphism. 
\ebew

\begin{rk}\em
Since $\msbar_{g,m}$ is a categorical 
quotient of $\Hgn$, there is a natural morphism of stacks 
$  \left[\Hgn/\PGL(N+1) \right]\rightarrow \msbar_{g,m}$. 
In fact, the diagram  
\[\diagram
\mgstack{g,m} \rrto^{\sim}\drto&& \left[\Hgn/\PGL(N+1) \right]\dlto \\
&  \msbar_{g,m} 
\enddiagram \]
commutes.
\end{rk}

\section{Moduli of $\mgamma$-pointed stable curves}

\begin{abschnitt}\em\label{gammaequiv}
In the introduction above we said that an $\mgamma$-pointed stable curve of genus $g$ over a closed point should be understood as an equivalence class of $m$-pointed stable curves $(f:C\rightarrow \Spec(k),\sigma_1,\ldots,\sigma_m)$, where $f:C\rightarrow \Spec(k)$ is an $m$-pointed stable curve with its marked points given by sections $\sigma_i: \Spec(k)\rightarrow C$ for $i=1,\ldots,m$, and where  the equivalence relation is taken with respect to permutations $\gamma\in \Gamma$ of the tuple $(\sigma_1,\ldots,\sigma_m)$. 

However, in order to provide a notion which is useful in applications, our definition \ref{mgdef} needs to be more technical. An $\mgamma$-pointed stable curve over an arbitrary scheme $S$ cannot simply be defined as an equivalence class of pointed stable curves $f:C\rightarrow S$, together with sections of marked points $\sigma_1,\ldots,\sigma_m: S \rightarrow C$. In fact, an $\mgamma$-pointed stable curve will in general not admit sections of marked points.

For example, choose a smooth stable curve $(f_0:C_0\rightarrow \Spec(k),P_1,\ldots,P_4)$ of genus $g=1$ with four marked points on it, for which there exists an automorphism $\tau:C_0\rightarrow C_0$, which interchanges the points $P_1$ and $P_2$, as well as $P_3$ and $P_4$. Let $U_1$ and $U_2$ be two Zariski open affine subschemes covering $\PP^1$. Obvioulsy the four distinguished points $P_1,\ldots,P_4$ on $C_0$ determine sections of the two direct products $p_1:U_1\times C_0\rightarrow U_1 $ and $p_2: U_2\times C_0\rightarrow U_2$. Thus the two products are both $4$-pointed stable curves of genus $1$ in the classical sense. We now glue the two trivial products with a twist, which is provided by applying the automorphism $\tau$ to the fibres over the intersection of $U_1$ and $U_2$. What we obtain is a family $f:C\rightarrow \PP^1$ of curves of genus $g=1$, each fibre with  $m=4$ distinguished points, which can not be labeled using global sections. In particular, this is no longer a $4$-pointed stable curve in the sense of Knudsen. However, there are two well-defined classes of distinguished points: those coming from points labeled $P_1$ or $P_2$ on the covering, and those coming from points $P_3$ or $P_4$. This is an example of an $\mgamma$-pointed stable curve of genus $g=1$ over $\PP^1$ as defined below, where $m=4$, and $\Gamma=\ZZ_2\times \ZZ_2$ is the subgroup of the permutation group $\Sigma_4$, which is generated by the two transpositions interchanging ``1'' and ``2'', and ``3'' and ``4'', respectively. The labels of the four distinguished points in each fibre are only determined up to permutations in $\Gamma$. 

In general, global sections of marked points exist only on an \'etale cover of the curve $f: C\rightarrow S$. Conversely, if there exist global sections on an \'etale cover, which satisfy a certain compatibility condition $(\ast)$, which will be discussed in remark \ref{class} in detail, then there will be well-defined classes of labels of the distinguished points.

In the case where $S = \Spec(k)$, lemma \ref{onetoone04} shows that the formal definiton agrees with our geometric intuition. 
\end{abschnitt}

\begin{defi}\em\label{mgdef}\label{n6}
Let $\Gamma$ be a subgroup of the symmetric group $\Sigma_m$. \\
$(i)$ A {\em charted $\mgamma$-pointed stable curve of genus $g$} is a tuple
\[ (f: C\rightarrow S, \:\: u: S' \rightarrow S, \:\: \sigma_1,\ldots,\sigma_m: S'\rightarrow C'),\]
where $f: C\rightarrow S$ is a flat and proper morphism of schemes, with reduced and connected algebraic curves as its fibres,  $u: S'\rightarrow S$ is an \'etale covering, defining  the fibre product $C' := C\times_S S'$, such that the induced morphism $f' : C' \rightarrow S'$, together with the sections $\sigma_1,\ldots,\sigma_m: S'\rightarrow C'$, is an $m$-pointed stable curve of genus $g$. Furthermore, we require that for all closed points $s,s'\in S'$ with $u(s)=u(s')$, there exists a permutation $\gamma_{s,s'}\in \Gamma$, such that 
\[ (\ast) \qquad \qquad \overline{u} \circ \sigma_i(s) = \overline{u}\circ \sigma_{\gamma_{s,s'}(i)}(s') \]
holds for all $i=1,\ldots,m$, where $\overline{u}: C'\rightarrow C$ is the morphism induced by $u: S' \rightarrow S$.\\
$(ii)$ Two charted $\mgamma$-pointed stable curves are called {\em equivalent}, if there exists a third charted $\mgamma$-pointed stable curve dominating both of them, in the sense of remark \ref{explain} below.\\
$(iii)$ An {\em $\mgamma$-pointed stable curve of genus $g$} is an equivalence class of  charted $\mgamma$-pointed stable curves of genus $g$.
\end{defi}

\begin{rk}\em\label{explain}
Let a charted $\mgamma$-pointed stable curve $(f: C\rightarrow S, \:\: u': S' \rightarrow S, \:\: \sigma_1,\ldots,\sigma_m: S'\rightarrow C')$ be given. A second charted $\mgamma$-pointed stable curve $(f: C\rightarrow S, \:\: u'': S'' \rightarrow S, \:\: \tau_1,\ldots,\tau_m: S''\rightarrow C'')$ with the same underlying curve $f:C\rightarrow S$ is said to {\em dominate} the first one, if the morphism $u'': S''\rightarrow S$ factors as $u'' = u' \circ v$, where $v: S''\rightarrow S'$ is an \'etale covering, which induces an isomorphism
\[ C'' \cong C'\times_{S'}S'', \]
and such that for all closed points $s\in S''$ there exists a permutation $\gamma_s\in \Gamma$, with
\[ (\ast') \qquad \qquad \overline{v}\circ \tau_i(s) = \sigma_{\gamma_s(i)}\circ v(s) \]
holding for all $i=1,\ldots, m$. Here $\overline{v}: C''\rightarrow C'$ denotes the morphism induced by $v: S'' \rightarrow S'$.
\end{rk}

\begin{rk}\em\label{class}
Let us analyse what the definition of an $\mgamma$-pointed stable curve means on the fibres of a curve, i.e. let us consider curves $f: C\rightarrow S$ for $S= \Spec(k)$.

$(i)$  Let $(f: C\rightarrow \Spec(k), \:\: u: S' \rightarrow \Spec(k), \:\: \sigma_1,\ldots,\sigma_m: S'\rightarrow C')$ be a charted  $\mgamma$-pointed stable curve of genus $g$ over a simple point. Then, by the compatibility condition $(\ast)$ of definition \ref{mgdef}, the sections  $\sigma_1,\ldots,\sigma_m$ distinguish $m$ distinct points on $C$, which are necessarily different from any nodes $C$ might have. 

Let $[i]$ denote the equivalence class of an element $i\in\{1,\ldots,m\}$ with respect to the action of $\Gamma$. By the compatibility condition $(\ast)$, the sections of $f': C'\rightarrow S'$ associate to each distinguished point of $C$ a unique label $[i]$, for some $i\in\{1,\ldots,m\}$. Formally we define for a distinguished point $q\in C$
\[ \mbox{class}(q) := [i], \]
if $q=\overline{u}\circ \sigma_i(s)$ for some $i\in \{1,\ldots,m\}$ and some $s\in S'$. In other words,  a charted  $\mgamma$-pointed stable curve of genus $g$ defines an $m$-pointed stable curve of genus $g$, but where the labels of the marked points are only determined up to a permutation in $\Gamma$. We will make this precise in lemma \ref{onetoone04} below. 

$(ii)$ Let $(f: C\rightarrow S, \:\: u'': S'' \rightarrow S, \:\: \tau_1,\ldots,\tau_m: S''\rightarrow C'')$ be a  second charted $\mgamma$-pointed stable curve with the same underlying curve $f:C\rightarrow S$, which dominates the first one. The compatibility condition $(\ast')$ of remark \ref{explain}, together with condition $(\ast)$ of definition \ref{mgdef}, ensures that the distinguished points on $C$ are the same in both cases. Even more, the classes of the distinguished points as defined above remain the same. 

$(iii)$ The notion of a charted $\mgamma$-pointed stable curve makes it necessary to specify one \'etale cover of $f:C\rightarrow S$ by an $m$-pointed stable curve. The definition of an $\mgamma$-pointed stable curve is independent of such a choice. 
\end{rk}

\begin{lemma}\label{onetoone04}
There is one-to-one correspondence between $\mgamma$-pointed stable curves of genus $g$ over $\Spec(k)$ and equivalence classes of $m$-pointed stable curves of the same genus over $\Spec(k)$ in the sense of remark \ref{gammaequiv}. 
\end{lemma}

\proof
Let $(f: C\rightarrow \Spec(k), \sigma_1,\ldots,\sigma_m)$ represent an equivalence class of $m$-pointed stable curves with respect to the action of $\Gamma$. Then $(f: C\rightarrow \Spec(k), \id_{\Spec(k)}, \sigma_1,\ldots,\sigma_m)$ is a charted  $\mgamma$-pointed stable curve. If $(f: C\rightarrow \Spec(k), \tau_1,\ldots,\tau_m)$ is a different representative of the above class, then by definition there exists a permutation $\gamma\in\Gamma$, such that $\sigma_i= \tau_{\gamma(i)}$ for all $i=1,\ldots,m$. The disjoint union $f \coprod f': C \coprod C' \rightarrow \Spec(k) \coprod \Spec(k)$ defines a charted $\mgamma$-pointed stable curve, dominating both $(f: C\rightarrow \Spec(k), \id_{\Spec(k)}, \sigma_1,\ldots,\sigma_m)$ and $(f: C\rightarrow \Spec(k), \id_{\Spec(k)}, \tau_1,\ldots,\tau_m)$. Hence to each equivalence class of $m$-pointed stable curves there is associated a well-defined $\mgamma$-pointed stable curve.   

Conversely, let an $\mgamma$-pointed stable curve be given by a representative 
 $(f: C\rightarrow \Spec(k), \:\: u: S' \rightarrow \Spec(k), \:\: \sigma_1,\ldots,\sigma_m: S'\rightarrow C')$. Choose some closed point $s\in S'$. Then the fibre $f_s': C_s'\rightarrow \Spec(k)$ of $f': C'\rightarrow S'$ over $s$, 
together with the sections $\tau_i := \sigma_i|\{s\}$, for $i=1,\ldots,m$,  defines an $m$-pointed stable curve $(f'_s: C_s'\rightarrow \Spec(k), \tau_1,\ldots,\tau_m)$. Let $(f'_{s'}: C_{s'}'\rightarrow \Spec(k), \tau_1',\ldots,\tau_m')$  denote the $m$-pointed curve obtained from the fibre over a different point $s'\in S'$. We clearly have $C_{s'}' \cong C \cong C_s'$. Because of condition $(\ast)$ of definition \ref{mgdef}  there exists a permutation $\gamma\in \Gamma$, such that $\tau_i= \tau_{\gamma(i)}$ for all $i=1,\ldots,m$. Hence the equivalence class of the $m$-pointed stable curve is well defined.

The two constructions are inverse to each other, which concludes the proof of  the lemma.
\ebew

\begin{rk}\em
Note that the permutations $\gamma_{s,s'} \in \Gamma$ in definition \ref{mgdef} are necessarily uniquely determined, and the same is true for the permutations $\gamma_s$ in remark \ref{explain}. 
\end{rk}

\begin{defi}\em 
$(i)$  Let $(f_k: C_k\rightarrow S_k, u_k: S_k'\rightarrow S_k, \sigma_1^{(k)},\ldots,\sigma_m^{(k)})$ be two charted $\mgamma$-pointed stable curves of genus $g$ for $k=1,2$. A {\em morphism between charted $\mgamma$-pointed stable curves} is a morphism of schemes 
\[ h : \:\: S_1 \rightarrow S_2,\]
which induces an isomorphism $C_1\cong C_2\times_{S_2} S_1$, and which satisfies the condition $(\diamond)$ below. 

Put $S_1'' := S_1' \times_{S_2} S_2'$. Then the natural composed morphism $u''_1: S_1''\rightarrow S_1$ is an \'etale covering of $S_1$, and for $C_1'':= C_1\times_{S_1} S_1''$ the tuple $(f_1: C_1\rightarrow S_1, u'':S_1''\rightarrow S_1, \overline{\sigma}_1^{(1)},\ldots,\overline{\sigma}_m^{(1)} : S_1''\rightarrow C_1'')$ is a charted $\mgamma$-pointed stable curve dominating $(f_1: C_1\rightarrow S_1, u':S_1'\rightarrow S_1, {\sigma}_1^{(1)},\ldots,{\sigma}_m^{(1)} : S_1'\rightarrow C_1')$. Here $\overline{\sigma}_i^{(1)}$ denotes the section induced by $\sigma_i^{(1)}$, for $i=1,\ldots,m$. There are induced morphisms $\overline{h}: S_1'' \rightarrow S_2'$ and $h'':C_1''\rightarrow C_2'$, such that the diagram
\[ \diagram
C_1'' \xto[rrr]^{h''} \xto[ddd] \drto &&& C_2' \xto[ddd] \dlto\\
& C_1 \rto \dto & C_2 \dto \\
& S_1 \rto^h & S_2 \\
S_1'' \xto[rrr]^{\overline{h}} \urto  &&& S_2' \ulto
\enddiagram \]
commutes, and all quadrangles are Cartesian. We require that there exists a permutation $\gamma\in \Gamma$, such that
\[ (\diamond) \qquad \qquad h''\circ \overline{\sigma}_i^{(1)} = \sigma_{\gamma(i)}^{(2)} \circ \overline{h} \]
holds for all $i=1,\ldots,m$. 

$(ii)$ A {\em morphism between $\mgamma$-pointed stable curves} is an equivalence class of morphisms between charted $\mgamma$-pointed stable curves.
\end{defi}

\begin{rk}\em
$(i)$ If $h: S_1 \rightarrow S_2$ is a morphism between a pair of charted $\mgamma$-pointed stable curves, then it is also a morphism for any other pair equivalent to it. This can easily be seen by drawing the appropriate extensions of the above commutative diagram. In particular, the notion of morphisms between $\mgamma$-pointed stable curves is well-defined.\\
$(ii)$ A morphism $h: S_1\rightarrow S_2$ between $\mgamma$-pointed stable curves is an isomorphism, if and only if it is an isomorphism of schemes $S_1\cong S_2$. Note however that not any isomorphism of schemes defines an isomorphism of $\mgamma$-pointed stable curves, it might not even define a morphism between them. \\
$(iii)$ Two $\mgamma$-pointed stable curves over the same underlying curve $f:C\rightarrow S$ are isomorphic if and only if they admit \'etale coverings by $m$-pointed stable curves, which are equivalent up to a permutation $\gamma$ of the labels of their sections, with $\gamma$ contained in  $\Gamma$, compare remark \ref{gammaequiv}. 
\end{rk}

\begin{rk}\em
A morphism $h: S_1\rightarrow S_2$ is a morphism of $\mgamma$-pointed stable curves if and only if there is a Cartesian diagram
\[ \diagram
C_1 \rto^{h'} \dto & C_2 \dto\\
S_1\rto^h &S_2
\enddiagram \]
such that for all closed points $s\in S_1$ and all distinguished points $q\in C_{1,s}$ in the fibre over $s$ holds
\[ \mbox{class}(q) = \mbox{class}(h'(q)) .\]
See remark \ref{class} for the notation. 
\end{rk}

\begin{rk}\em
Recall that the symmetric group $\Sigma_m$ acts naturally on the set of $m$-pointed stable curves by permutation of the labels of the marked points. Since the $m$ sections of marked points are disjoint, the action of $\Sigma_m$, and of any subgroup $\Gamma \subset \Sigma_m$,  on $\overline{H}_{g,n,m}$ is free. 
\end{rk}

\begin{rk}\em
A subgroup $\Gamma$ of $\Sigma_m$ acts on $\msbar_{g,m}$. The quotient 
\[ \msbar_{g,\mgamma} := \msbar_{g,m} / \Gamma \label{n7}\]
is a coarse moduli space for $\mgamma$-pointed curves.  Let $f: C\rightarrow S$ be an $\mgamma$-pointed stable curve of genus $g$. By lemma \ref{onetoone04}, the closed points of $\msbar_{g,\mgamma}$ are in one-to-one correspondence with isomorphism classes of $ \mgamma$-pointed stable curves of genus $g$ over $\Spec(k)$.

Let $(f:C\rightarrow S, u:S'\rightarrow S, \sigma_1,\ldots,\sigma_m: S'\rightarrow C')$ be a charted $\mgamma$-pointed stable curve of genus $g$. Since $f': C'\rightarrow S'$ is an $m$-pointed stable curve, there is an induced morphism $\theta_{f'}:S'\rightarrow  \msbar_{g,m}$, and by composition a morphism $\overline{\theta}_{f'}:S'\rightarrow  \msbar_{g,\mgamma}$. By definition, for any closed point $s\in S$, and for any pair of points $s',s''\in S''$ with $u(s')=u(s'')=s$, holds the identity $\overline{\theta}_{f'}(s') =\overline{\theta}_{f'}(s'')$. Thus there is an induced morphism $\theta: S\rightarrow \msbar_{g,\mgamma}$, such that $\theta\circ u =  \overline{\theta}_{f'}$.  Clearly any other charted $\mgamma$-pointed stable curve dominating this one induces the same morphism from $S$ into $\msbar_{g,\mgamma}$. Hence there is a well defined morphism $\theta_f: S\rightarrow \msbar_{g,\mgamma}$ for any $\mgamma$-pointed stable curve $f: C\rightarrow S$. Furthermore, by referring back to $\msbar_{g,m}$ again, one easily sees that $\msbar_{g,\mgamma}$ dominates any other scheme with this universal property, so it is indeed a coarse moduli space for $\mgamma$-pointed stable curves.
\end{rk}

\begin{rk}\em
If $\Gamma = \{\id\}$ is the trivial subgroup of of $\Sigma_m$, then the glueing condition $(\ast)$ for $\mgamma$-pointed stable curves implies that there are  $m$ induced  global sections of marked points of $f:C\rightarrow S$.  Thus in this case an $\mgamma$-pointed stable curve is just an $m$-pointed stable curve in the classical sense.
\end{rk}

\begin{rk}\em\label{rk224}
Let $\Gamma \subset \Sigma_m$ be a subgroup. There is a free action of $\Gamma$ on $\Hgn$. The quotient shall be denoted by
\[ \Hgmp := \Hgn / \Gamma \label{n8}.\]
Let $u_{g,n,m} : {\cal C}_{g,n,m} \rightarrow \Hgn$ denote the universal curve over $\Hgn$. Note that fibres of $u_{g,n,,m}$ over such points of $\Hgn$, which correspond to each other under the action of $\Gamma$, are identical as embedded curves in $\PP^N$, with the same distinguished points on them. The only difference is in the labels of the marked points, which are permuted by elements of $\Gamma$. Thus the action of $\Gamma$ on $\Hgn$ extends to an action on ${\cal C}_{g,n,\mgamma}$, by identifying fibres over corresponding points. 
The quotient
\[ u_{g,n,\mgamma}: \quad {\cal C}_{g,n,\mgamma} := {\cal C}_{g,n,m}/ \Gamma \: \longrightarrow \: \Hgn/\Gamma = \Hgmp\]
exists as an $\mgamma$-pointed stable curve over $\Hgmp$. The original universal curve $u_{g,n,m} : {\cal C}_{g,n,m} \rightarrow \Hgn$, together with the \'etale morphism $\Hgn\rightarrow \Hgmp$ defines a charted $\mgamma$-pointed stable curve representing it. 
\end{rk}

\begin{abschnitt}\em
Let $f : C\rightarrow S$ be an $\mgamma$-pointed stable curve of genus $g$, which is given together with an embedding of the underlying curve $f: C\rightarrow S$ into $\pr_2: \PP^N\times S \rightarrow  S$. For any \'etale covering $u:S'\rightarrow S$, there is an induced embedding of $f' : C' = C\times_S S' \rightarrow S'$ into $\pr_2: \PP^N\times S' \rightarrow S'$. Hence, for a representing charted $\mgamma$-pointed curve, the curve $f' : C' \rightarrow S'$ is an embedded $m$-pointed stable curve, and thus induces a morphism 
\[ \theta_{f'}: \: \: S' \rightarrow \Hgn.\]
 For any closed point $s\in S$, the embedding of the fibre $C'_{s'}$ into $\PP^N$ is the same for  all points $s'\in S'$ with $u(s')=s$, as it is nothing else but the embedding of $C_s$ into $\PP^N$. Note that the points in $\Hgn$ representing different fibres $C'_{s'}$ need not be the same, as the labels of the marked points may differ by a permutation in $\Gamma$. However, the morphism $\overline{\theta}_{f'}: S' \rightarrow \Hgmp$ obtained by composition with the quotient map, factors through a morphism $\theta_f : S\rightarrow \Hgmp$. This morphism is even independent of the chosen charted $\mgamma$-pointed stable curve representing $f:C\rightarrow S$. 

By the universal property of the Hilbert scheme, there is an isomorphism of $m$-pointed stable curves between $\theta_{f'}^\ast {\cal C}_{g,n,m}$ and $C'$ over $S'$. As a scheme over $S'$, the curve $\theta_{f'}^\ast {\cal C}_{g,n,m}$ is isomorphic to  $\overline{\theta}_{f'}^\ast {\cal C}_{g,n,\mgamma}$, even though the latter has no natural structure as an $m$-pointed curve. Therefore, by the definition of $C' = u^\ast C$, and since $\overline{\theta}_{f'}= \theta_f\circ u$, one can construct an isomorphism of schemes over $S$ between $C$ and $\theta_f^\ast {\cal C}_{g,n,\mgamma}$. Together with the isomorphisms between the covering $m$-pointed stable curves, this shows that $f:C\rightarrow S$ is isomorphic to $\theta_f^\ast {\cal C}_{g,n,\mgamma} \rightarrow S$ as an $\mgamma$-pointed stable curve. 

In other words,  $\Hgmp$ is in fact a
fine moduli space for embedded $\mgamma$-pointed stable curves of genus $g$, with universal curve $u_{g,n,\mgamma}: {\cal C}_{g,n,\mgamma} \rightarrow \Hgmp$. 
\end{abschnitt}

\begin{rk}\em
The action of $\PGL(N+1)$ on $\Hgn$ induces an action of $\PGL(N+1)$ on $\Hgmp$. Indeed, by the construction of the universal embedded curve over  $\Hgn$, the embedding of one of its fibres, considered as an $m$-pointed stable curve of genus $g$, does not depend on the ordering of the labels of its marked points. So the action of $\PGL(N+1)$ on $\Hgn$ commutes with the action of any subgroup $\Gamma\subset \Sigma_m$ on $\Hgn$.  The coarse moduli space $\msbar_{g,\mgamma}$ is a GIT-quotient of $\Hgmp$ by the action of $\PGL(N+1)$. 
\end{rk}

\begin{defi}\em
The {\em moduli stack $\, \mgstack{g,\mgamma}$\label{n9} of $\mgamma$-pointed stable curves of genus $g$} is the stack defined as the category fibred in groupoids over  the category of schemes, where for a scheme $S$ the objects in the fibre category $\mgstack{g,\mgamma}(S)$ are the $\mgamma$-pointed stable curves of genus $g$ over $S$.
Morphisms in $\mgstack{g,\mgamma}$ are morphisms between $\mgamma$-pointed stable curves.
\end{defi}

\begin{prop}\label{P217}
There is an isomorphism of stacks
\[  \mgstack{g,\mgamma} \: \cong \: \left[ \, \Hgmp / \PGL(N+1) \right] .\]
\end{prop}

\begin{rk}\em
Note that by the general theory of quotient stacks, and the fact that the action of $\Gamma$ on $\Hgn$ is free, there is an isomorpism of stacks
\[   \left[ \, \Hgmp / \PGL(N+1) \, \right]  \: \cong \: \left[ \, \Hgn / \Gamma \times \PGL(N+1) \right]
 .\]
In particular, the quotient stack  $  \mgstack{g,\mgamma}$ can be viewed as a stack quotient of the moduli stack $  \mgstack{g,m}$ with respect to the group $\Gamma$. 
\end{rk}

\proofof{the proposition}
The proof is analogous to that of proposition \ref{28}, so we may be brief here, and concentrate on those parts of the proof, which need to be adapted. \\
$(i)$ Let $f:C\rightarrow S \in \ob(\mgstack{g,\mgamma})$ be an $\mgamma$-pointed stable curve. For each 
\'etale cover $u: S'\rightarrow S$, the $m$ disjoint sections $\sigma_1,\ldots,\sigma_m: S'\rightarrow C'$ define hypersurfaces $S_1,\ldots,S_m$ in $C'$. Their images in $C'$ define a divisor in $C$, which will be denoted by $S_{\mgamma}$. Note that this divisor is independent of the chosen \'etale cover. 

 Let $p:E\rightarrow S$ be the principal $\PGL(N+1)$-bundle associated to the projective bundle $\PP f_\ast( \omega_{C/S}(S_{\mgamma}))^{\otimes n}$. Let $E' := E\times_S S'$ be the pullback of $E$ to some \'etale covering $u:S'\rightarrow S$ of $S$. In the same way as in the proof of proposition \ref{28}, there is a natural embedding of the $m$-pointed stable curve $C'\times_{S'} E' \rightarrow E'$, and this induces a morphism
\[ \phi' : \: E' \rightarrow \Hgn .\]
By composition with the quotient map there is a morphism
\[\overline{\phi} : \:\: E' \rightarrow \Hgn/\Gamma = \Hgmp , \]
which is $\PGL(N+1)$-equivariant. Note that the action of $\Gamma$ commutes with the action of $\PGL(N+1)$ on $\Hgn$.  

For all closed points $e,e' \in E'$ which project to the same point in $E$, we have $\overline{\phi}(e)=\overline{\phi}(e')$. Therefore  $\overline{\phi}$ factors through a morphism $\phi: E\rightarrow \Hgmp$, which is also $\PGL(N+1)$-equivariant. Thus we obtain an object $(E,p,\phi)\in [\Hgmp/\PGL(N+1](S)$. 

The construction of the functor from $\mgstack{g,\mgamma}$ to $[\Hgmp/\PGL(N+1]$ on morphisms is straightforward. 

$(ii)$ Conversely, consider a triple  $(E,p,\phi)\in [\Hgmp/\PGL(N+1](S)$ for some scheme $S$. The morphism $\phi: E \rightarrow \Hgmp$ determines an embedded $\mgamma$-pointed stable curve $f': C' \rightarrow E$ of genus $g$, together with an isomorphism between $C'$ and the pullback of the universal curve ${\cal C}_{g,n,\mgamma}$ as schemes over $E$. The quotient 
\[ f: C := C' / \PGL(N+1) \: \rightarrow  \: E / \PGL(N+1) \cong S \]
exists, and it is by construction an $\mgamma$-pointed stable curve. 

One verifies that the two functors, whose  construction we outlined  above,  are inverse to each other as morphisms of stacks.
\ebew

\begin{rk}\em
From proposition \ref{P217} it follows in particular that $\mgstack{g,\mgamma}$ is a smooth Deligne-Mumford stack, with $\msbar_{g,\mgamma}$ as its moduli space.
 There is a canonical commutative diagram
\[ \diagram
\mgstack{g,m} \dto^\cong \rto &\mgstack{g,\mgamma} \dto^\cong\\
\left[\Hgn/\PGL(N+1) \right]\dto \rto & \left[\Hgn/\Gamma\times \PGL(N+1) \right]\dto \\
 \msbar_{g,m} \rto & \msbar_{g,\mgamma}.
\enddiagram \]

All horizontal arrows represent finite morphisms of degree equal to the order of $\Gamma$. The morphisms between the stacks are unramified, while in general the morphisms between the moduli spaces are not.  
\end{rk}

\begin{rk}\em
Note that for everything we said above analogous statements hold true if we replace $\Hgn$ 
and $\msbar_{g,m}$ by the open subschemes $H_{g,n,m}$ and $\ms_{g,m}$, as well as  $\Hgmp$ 
and $\msbar_{g,\mgamma}$ by the open subschemes $H_{g,n,\mgamma}$ and 
$\ms_{g,\mgamma}$, respectively.
\end{rk}

\end{document}